\documentclass[12pt,reqno]{amsart}
\usepackage{amssymb,amsthm,url}
\usepackage[breaklinks=true]{hyperref}
\usepackage[utf8]{inputenc}
\usepackage[OT2,T1]{fontenc}
\swapnumbers
\let\<\langle
\let\>\rangle

\usepackage{lscape}
\usepackage{tikz-cd}

\DeclareSymbolFont{cyrletters}{OT2}{wncyr}{m}{n}
\DeclareMathSymbol{\Sha}{\mathalpha}{cyrletters}{"58}

\theoremstyle{definition}
\newtheorem{definition}{Definition}[section]
\newtheorem{remark}[definition]{Remark}
\newtheorem{example}[definition]{Example}

\theoremstyle{plain}
\newtheorem{lemma}[definition]{Lemma}
\newtheorem{proposition}[definition]{Proposition}
\newtheorem{corollary}[definition]{Corollary}
\newtheorem{conjecture}[definition]{Conjecture}
\newtheorem{theorem}[definition]{Theorem}

\newenvironment{Proof}[1][Proof.]{\begin{trivlist}
\item[\hskip \labelsep {\bfseries #1}]}{\flushright
$\Box$\end{trivlist}}

\usepackage[all]{xy}
\usepackage{pdfsync}
\usepackage{ytableau}

\newcommand{\sudda}[1]{}

\begin{document}

\title{On the transposed Poisson $n$-Lie algebras}

\author{Farukh Mashurov}

\address{SDU University, Kaskelen, Kazakhstan.}
\address{Shenzhen International Center for Mathematics, Southern University of Science and Technology, Shenzhen, 518055, China.}

\email{f.mashurov@gmail.com}

\keywords{$n$-Lie algebra, transposed Poisson algebra, simple algebra, polynomial identity}
\subjclass[2020]{15A15, 17A30, 17A42, 17B63}


\maketitle

\begin{abstract} We study unital commutative associative algebras and their associated $n$-Lie algebras, showing that they are strong transposed Poisson $n$-Lie algebras under specific compatibility conditions. Furthermore, we generalize the simplicity criterion for transposed Poisson algebras, proving that a transposed Poisson $n$-Lie algebra is simple if and only if its associated $n$-Lie algebra is simple. In addition, we study the strong condition for transposed Poisson $n$-Lie algebras, proving that it fails in the case of a free transposed Poisson $3$-Lie algebra.
\end{abstract}

\section{\label{nn}\ Introduction}

In 1985, V.T. Filippov introduced $n$-Lie algebras as a generalization of the classical Lie algebras \cite{Fil1985}. These algebras have attracted attention due to their close connection with Nambu mechanics \cite{Nambu, Takhtajan}. A generalization of classical Hamiltonian mechanics was proposed by Nambu in 1973 \cite{Nambu}.

Filippov gave two important examples of $n$-Lie algebras \cite{Fil1985, Fil1998}. The first example is the vector product algebra  $A_{n+1}$, which is an  $(n + 1)$-dimensional $n$-ary algebra. The second example is the Jacobian algebra, which is a commutative associative algebra $A$ that with fixed pairwise commuting derivations $D_1, \ldots, D_n.$ In this algebra, the multiplication is defined by the determinant of the Jacobian matrix (see Example \ref{ex: Jac}).

Dzhumadil'daev provided another example  (Example \ref{ex: W}) of an infinite-dimensional $n$-Lie algebra in \cite{DzhJac}.  All simple linearly compact $n$-Lie superalgebras and (generalized) $n$-Nambu-Poisson algebras over a field $F$ of characteristic zero were classified by Cantarini and Kac \cite{CanKac2010, CanKac2016}. They proved that every simple linearly compact generalized $n$-Nambu-Poisson algebra is either isomorphic to the one presented by Filippov or Dzhumadil'daev algebra \cite{CanKac2016}.


A transposed Poisson algebra (or transposed Poisson $n$-Lie algebra) was recently introduced in \cite{BBGW2023} as a dualisation of the Poisson ($n$-Nambu-Poisson) algebra. It is defined by reversing the two binary operation roles in the Leibniz rule for a classical Poisson algebra.


The theory of transposed Poisson ($n$-Lie) algebras has attracted more attention in recent studies; see \cite{BBGW2023, BFK2024, CB2024, DIM25, Ouar2023, Sar2024} and the references therein. In particular, the authors in \cite{BFK2024} provided an overview of the known results in this area and a list of open questions.
One of these questions was about the examples and classifying simple transposed Poisson $n$-Lie algebras (see  \cite[Question 38]{BFK2024}). In \cite{Ouar2023}, it was proved that a transposed Poisson ($2$-Lie)  algebra $(L,\cdot,[\cdot,\cdot])$ is simple if and only if its associated Lie algebra $(L,[\cdot,\cdot])$ is simple.

 Dzhumadil'daev introduced the Poisson $n$-Lie algebra, which can be used to construct the $(n + 1)$-Lie algebra under an additional condition, called the strong condition \cite{DzhJac}.
The authors in \cite{BBGW2023} proved that this strong condition holds for the transposed Poisson 2-Lie algebra. They also constructed 3-Lie algebras from transposed Poisson 2-Lie algebras using derivations. This construction motivated them  to formulate the following conjecture:
\begin{conjecture}\label{conj}
    Let $n\geq 2$ be an integer. Let $(L,\cdot,[\cdot,\ldots,\cdot])$ be a transposed Poisson $n$-Lie algebra and let $D$ be a derivation of $(L,\cdot)$ and $(L,[\cdot,\ldots,\cdot]).$  Define an $(n + 1)$-ary operation
    $$\mu(a_1,\ldots,a_{n+1})=\sum_{i=1}^{n+1}(-1)^{i-1}D(a_i)[a_1,\ldots,\hat{a}_i,\ldots,a_{n+1}],$$
    where $\hat{a}_i$ means that the $i$-th entry is omitted. Then $(L,\cdot,\mu(\cdot,\ldots,\cdot))$  is a transposed Poisson $(n+1)$-Lie algebra.
\end{conjecture}

In \cite{HCCD}, the authors proved that Conjecture \ref{conj} holds for strong transposed Poisson $n$-Lie algebras. Moreover, a classification of transposed Poisson 3-Lie algebras of dimension 3 is given in \cite{YCL25}.

First, we study a commutative associative algebra $(A, \cdot)$ and an $n$-Lie algebra $(A, [\cdot,\ldots,\cdot])$ related by specific compatibility conditions.   Under these conditions, we prove that the algebra $(A,\cdot, [\cdot,\ldots,\cdot])$ is a strong transposed Poisson $n$-Lie algebra. As a corollary, we establish that the $n$-Lie algebra introduced by  Dzhumadil'daev with associative commutative multiplication is a strong transposed Poisson $n$-Lie algebra. Next, we generalize the result of  \cite{Ouar2023} to the $n$-ary case and prove that a transposed Poisson $n$-Lie algebra $(L,\cdot,[\cdot,\dots,\cdot])$ is simple if and only if the associated $n$-Lie algebra $(L,[\cdot,\dots,\cdot])$ is simple. Using our results combined with the classification of \cite{CanKac2016}, we provide an example of a simple linearly compact transposed Poisson $n$-Lie algebra. 
The authors of \cite{HCCD} noted that there is no evidence disproving the validity of the strong condition for all transposed Poisson $n$-Lie algebras with $n \geq 3$. In the last section, we show that the strong condition fails in the case $n = 3$ for the free transposed Poisson $n$-Lie algebra.

We work over the field of complex numbers $\mathbb{C}$, though all results hold over any field
of characteristic zero.

 \section{\label{preliminaries}\ Preliminaries}

In this section, we first recall the necessary definitions and some known results.
\begin{definition}
 An $n$-Lie algebra over a field $F$ is a vector space $L$ equipped with an $n$-ary skew-symmetric bracket operation:
\[
[\cdot,  \dots, \cdot] : L \times  \dots \times L\rightarrow L,
\]
such that for all $x_1, \dots, x_{n}, y_2, \dots, y_n \in L$, the following \emph{generalized Jacobi identity} holds:
\begin{equation}\label{id: jacobi}
    [[x_1,  \dots, x_{n}],y_2, \dots, y_{n}]=\sum_{i=1}^n [x_1,  \dots, x_{i-1},[x_i,y_2, \dots, y_{n}],x_{i+1},\ldots,x_n].\end{equation}
\end{definition}

There are two well-known infinite-dimensional $n$-Lie algebras, one of which  introduced by Filippov \cite{Fil1985, Fil1998}:

\begin{example}\label{ex: Jac}
 Let $A$ be a commutative associative algebra  with fixed pairwise commuting
derivations $\partial_1,\ldots,\partial_n$ of the algebra $A.$ Define an $n$-ary multiplication as follows:
$$Jac(u_1,u_2,\ldots,u_n)=\begin{array}{|ccc|}
  \partial_1 u_1   &  \cdots & \partial_1 u_n \\
    \cdots & \cdots & \cdots \\

      \partial_n u_1   &  \cdots & \partial_n u_n
\end{array},$$
for every $u_1,u_2,\ldots,u_n\in A.$
Then $(A, Jac(\cdot,\ldots,\cdot))$ is an $n$-Lie algebra. 

Moreover, this algebra satisfies the following Leibniz rule:
\begin{equation}\label{id: Jac leib} Jac(u_1 u_2,u_3,\ldots,u_{n+1})=u_1Jac( u_2,u_3,\ldots,u_{n+1})+u_2Jac(u_1 ,u_3,\ldots,u_{n+1}).\end{equation}

\end{example}

Another important example of an infinite-dimensional $n$-Lie algebra was given by Dzhumadil'daev in \cite{DzhJac}:
\begin{example}\label{ex: W}
Let $A$ be a commutative associative algebra $A$ with  $n-1$  pairwise commuting derivations $D_1,\ldots,D_{n-1}$. Define an $n$-ary bracket by
$$W(u_1,u_2,\ldots,u_n)=\begin{array}{|ccc|}
    u_1   &  \cdots &  u_n \\
  D_1 u_1   &  \cdots & D_1 u_n \\
    \cdots & \cdots & \cdots \\
     D_{n-1} u_1   &  \cdots & D_{n-1} u_n
\end{array},$$ 
for every $u_1,u_2,\ldots,u_n\in A.$ Then $(A, W(\cdot,\ldots,\cdot))$ is an $n$-Lie algebra  \cite{DzhJac}. 
\end{example}

Next, we introduce the notion of a \emph{transposed Poisson $n$-Lie algebra}, as defined in \cite{BBGW2023}.

\begin{definition}
    An algebra $(A,\cdot, [\cdot,\ldots,\cdot])$ is called transposed Poisson $n$-Lie algebra if 
    \begin{itemize}
        \item $(A,\cdot)$  is a  commutative associative algebra.
        \item $(A,[\cdot,\ldots,\cdot])$  is an $n$-Lie algebra.
        \item  for arbitrary elements $h, a_1,\ldots, a_n \in A$ the following Leibniz rule holds: \begin{equation}\label{id: tr P} n \, h\, [a_1,\ldots,a_n]=\sum_{i=1}^{n} [a_1,\ldots h a_i,\ldots, a_n]. \end{equation}  
    \end{itemize}
\end{definition}

When $n = 2,$ the transposed Poisson 2-Lie algebras are called transposed Poisson algebras, as introduced in \cite{BBGW2023}.

The following important identities hold in transposed Poisson $n$-Lie algebras and were established in \cite{HCCD}.

\begin{proposition} The following identities hold true in a transposed Poisson $n$-Lie algebra $(A,\cdot, [\cdot,\ldots,\cdot]):$ 

    \begin{equation}\label{id: [[]h]}
        [[x_1, \dots, x_{n}]h,y_2, \dots, y_{n}]=\sum_{i=1}^n (-1)^{i-1} [[x_i,y_2, \dots, y_{n}]h,x_1, \dots, \hat{x_{i}}, \ldots,x_n],
    \end{equation}
    \begin{equation}\label{id: h[)}
    \sum^{n+1}_{i=1}(-1)^{i-1}x_i [x_1,\ldots,\hat{x_i},\ldots,x_{n+1}] = 0,
    \end{equation}  where $h,x_i,y_i\in  A$ and $i=1,\ldots,n+1.$
    \end{proposition}

Following \cite{HCCD}, we now introduce the definition of a \emph{strong} transposed Poisson $n$-Lie algebra.

\begin{definition} A transposed Poisson $n$-Lie algebra $(A, \cdot, [\cdot, \ldots , \cdot])$ is called strong if the following identity holds:
   $$y_1 [h y_2,x_1,\ldots ,x_{n-1}]-y_2 [h y_1,x_1,\ldots,x_{n-1}]\\
    +\sum_{i=1}^{n-1}(-1)^{i-1}h x_i[y_1,y_2,x_1,\ldots,\hat{x}_i,\ldots,x_{n-1}]=0,$$
for any $h,y_1,y_2,x_i\in A,1\leq i\leq n-1.$
\end{definition}

Note that when $n = 2,$ the above identity has the following form
$$y_1 [hy_2, x_1] + y_2 [x_1, hy_1] + hx_1 [y_1, y_2] = 0.$$
In this case, the strong condition holds true in every transposed Poisson (2-Lie) algebra by Theorem 2.5 (11) in \cite{BBGW2023}.

\section{Unital transposed Poisson $n$-Lie algebras}

 In this section, we study a unital commutative associative algebra $(A, \cdot)$  and an $n$-Lie algebra $(A, [\cdot,\ldots,\cdot])$  defined by specific compatibility conditions.   We prove that this algebra $(A,\cdot, [\cdot,\ldots,\cdot])$ is a strong transposed Poisson $n$-Lie algebra.

The following proposition is a direct consequence of the symmetry of the product $ \cdot$ and the skew-symmetry of the bracket $ [ \cdot, \ldots, \cdot] $.

\begin{proposition}  Let $(A, \cdot)$ be a unital  commutative algebra and  $(A, [\cdot,\ldots,\cdot])$ be  an $n$-ary skew-symmetric bracket.
Then the following identities hold true in $(A,\cdot, [\cdot,\ldots,\cdot]):$ 
\begin{equation}\label{id: a_ia_j[]}
\sum_{i=1}^{n}\sum_{j=1,j\neq i}^n (-1)^{j-1} a_{i} a_j[1,a_1,\ldots,\hat{a_j} ,\ldots, a_{n+1}\hat{a_i},\ldots, a_n]=0,\end{equation}

\begin{equation}\label{id: ha_ia_j[]}
\sum_{i=1}^{n}(-1)^{i-1}\sum_{j=1,j\neq i}^n (-1)^{j-1} a_{i} a_j[1,a_{n+1},a_1,\ldots,\hat{a_j} ,\ldots, \hat{a_i},\ldots, a_n]=0,\end{equation}
where $a_1,\ldots,a_n\in  A,$ , where $\hat{a_l}$ indicates that element $a_l$ is omitted, and $a_{n+1}\hat{a_k}$ means that $a_{n+1}$ is at position $k$.  
\end{proposition}
\begin{Proof}
To prove the first identity \eqref{id: a_ia_j[]}, consider any pair of indices $i > j$. By the skew-symmetry of the $n$-ary bracket, we observe the following relation:  
\begin{align*}
    (-1)^{i-1}a_{j} a_i[1,a_1,\ldots,a_{n+1}\hat{a_j} ,\ldots, \hat{a_i},\ldots, a_n] & = (-1)^{i-1+(i-j-1)}a_{j} a_i[1,a_1,\ldots,\hat{a_j} ,\ldots, a_{n+1}\hat{a_i},\ldots, a_n]\\
    & = (-1)^{j}a_{j} a_i[1,a_1,\ldots,\hat{a_j} ,\ldots, a_{n+1}\hat{a_i},\ldots, a_n].
\end{align*}

In the expansion of the sum in \eqref{id: a_ia_j[]}, we encounter the following paired terms:
$$(-1)^{j-1}a_{i} a_j[1,a_1,\ldots,\hat{a_j} ,\ldots, a_{n+1}\hat{a_i},\ldots, a_n]+(-1)^{i-1}a_{j} a_i[1,a_1,\ldots,a_{n+1}\hat{a_j} ,\ldots, \hat{a_i},\ldots, a_n].$$

On the other hand, applying the previously established relation to the second term, we obtain:
$$(-1)^{j-1}a_{i} a_j[1,a_1,\ldots,\hat{a_j} ,\ldots, a_{n+1}\hat{a_i},\ldots, a_n]+(-1)^{i-1}a_{j} a_i[1,a_1,\ldots,a_{n+1}\hat{a_j} ,\ldots, \hat{a_i},\ldots, a_n]$$
$$(-1)^{j-1}a_{i} a_j[1,a_1,\ldots,\hat{a_j} ,\ldots, a_{n+1}\hat{a_i},\ldots, a_n]+(-1)^{j}a_{j} a_i[1,a_1,\ldots,\hat{a_j} ,\ldots, a_{n+1}\hat{a_i},\ldots, a_n].$$

Since $(-1)^{j-1} + (-1)^{j} = 0$, the above expression is zero, implying that each term cancels in the sum. Therefore, we conclude that  
$$\sum_{i=1}^{n}\sum_{j=1,j\neq i}^n (-1)^{j-1} a_{i} a_j[1,a_1,\ldots,\hat{a_j} ,\ldots, a_{n+1}\hat{a_i},\ldots, a_n]=0.$$
An analogous argument establishes the second identity \eqref{id: ha_ia_j[]}, completing the proof. 
\end{Proof}


We now prove the main theorem of this section: 

\begin{theorem}
    \label{th: if bracket with 1 the Tr}
    Let $(A, \cdot)$ be a unital  commutative algebra and  $(A, [\cdot,\ldots,\cdot])$ be  an $n$-Lie algebra such that
for every $u_1,\ldots,u_n\in A$ we have

\begin{equation}\label{br as lcomb 1}
[u_1,\ldots,u_n]=\sum^{n}_{i=1}(-1)^{i-1}u_i [1,u_1,\ldots,\hat{u_i},\ldots,u_{n}],\end{equation}
and 
\begin{equation}\label{leib w 1}
[1,u_1u_2,u_3\ldots,u_n]=u_1[1,u_2,u_3\ldots,u_n]+u_2[1,u_1,u_3\ldots,u_n].\end{equation}
Then $(A,\cdot, [\cdot,\ldots,\cdot])$  is a strong transposed Poisson $n$-Lie algebra. 
\end{theorem}
\begin{Proof} We first show that $(A,\cdot, [\cdot,\ldots,\cdot])$ with identities \eqref{br as lcomb 1} and \eqref{leib w 1} is a transposed Poisson $n$-Lie algebra. That is, the identity $$n \, a_{n+1}\, [a_1,\ldots,a_n]=\sum_{i=1}^{n} [a_1,\ldots ,a_{n+1} a_i,\ldots, a_n]$$ holds true in $(A,\cdot, [\cdot,\ldots,\cdot]),$ where $a_1,\ldots,a_{n+1}\in A.$ 

 Applying \eqref{br as lcomb 1} to the sum  $\sum_{i=1}^{n} [a_1,\ldots, a_{n+1} a_i,\ldots, a_n],$ we have
$$\sum_{i=1}^{n} [a_1,\ldots a_{n+1} a_i,\ldots, a_n]=$$
$$\sum_{i=1}^{n}(-1)^{i-1}a_{n+1} a_i[1,a_1,\ldots,\hat{a_i},\ldots, a_n]+\sum_{i=1}^{n}\sum_{j=1,j\neq i}^n (-1)^{j-1} a_j[1,a_1,\ldots,\hat{a_j},\ldots, a_{n+1} a_i,\ldots, a_n]=$$
(applying \eqref{leib w 1} to the second sum we have)
$$ \sum_{i=1}^{n}(-1)^{i-1}a_{n+1} a_i[1,a_1,\ldots,\hat{a_i},\ldots, a_n]+\sum_{i=1}^{n}\sum_{j=1,j\neq i}^n (-1)^{j-1} a_{n+1} a_j[1,a_1,\ldots,\hat{a_j},\ldots, a_i,\ldots, a_n]+$$$$\sum_{i=1}^{n}\sum_{j=1,j\neq i}^n (-1)^{j-1} a_{i} a_j[1,a_1,\ldots,\hat{a_j},\ldots,  a_{n+1} \hat{a_i},\ldots, a_n].$$

Note that, by \eqref{id: a_ia_j[]} the last sum is equal to zero in the above expression. Thus, we have 
$$\sum_{i=1}^{n} [a_1,\ldots a_{n+1} a_i,\ldots, a_n]=$$
$$ \sum_{i=1}^{n}(-1)^{i-1}a_{n+1} a_i[1,a_1,\ldots,\hat{a_i},\ldots, a_n]+\sum_{i=1}^{n}\sum_{j=1,j\neq i}^n (-1)^{j-1} a_{n+1} a_j[1,a_1,\ldots,\hat{a_j},\ldots, a_i,\ldots, a_n].$$
Moreover, the expression above can be rewritten as follows:
\begin{align*}
      \sum_{i=1}^{n} [a_1,\ldots a_{n+1} a_i,\ldots, a_n]& =\sum_{i=1}^{n}\sum_{j=1}^n (-1)^{j-1} a_{n+1} a_j[1,a_1,\ldots,\hat{a_j},\ldots, a_n]\\
     & =n \; a_{n+1} \big (\sum_{j=1}^n (-1)^{j-1}  a_j[1,a_1,\ldots,\hat{a_j},\ldots, a_n] \big )\\
     &\stackrel{\mathrm{\eqref{br as lcomb 1}}}{=} n \, a_{n+1}\, [a_1,\ldots,a_n].
\end{align*}

Therefore, we have
    $$\sum_{i=1}^{n} [a_1,\ldots ,a_{n+1} a_i,\ldots, a_n] = n \, a_{n+1}\, [a_1,\ldots,a_n].$$

Hence, $(A,\cdot, [\cdot,\ldots,\cdot])$  is a transposed Poisson $n$-Lie algebra. 

Next, we show that  this  transposed Poisson $n$-Lie algebra $(A,\cdot, [\cdot,\ldots,\cdot])$ is a strong. That is the following identity holds true    $$y_1 [h y_2,a_1,\ldots ,a_{n-1}]-y_2 [h y_1,a_1,\ldots,a_{n-1}]\\
    +\sum_{i=1}^{n-1}(-1)^{i-1}h a_i[y_1,y_2,a_1,\ldots,\hat{a}_i,\ldots,a_{n-1}]=0,$$ where $h,y_1,y_2, a_1,\ldots,a_{n+1}\in A.$ 
    
    Applying \eqref{br as lcomb 1} and \eqref{leib w 1} to $y_1 [h y_2,a_1,\ldots ,a_{n-1}]$  we have  
 \begin{align*}
     y_1 [h y_2,a_1,\ldots ,a_{n-1}]& \stackrel{\mathrm{\eqref{br as lcomb 1}}}{=} y_1 h y_2 [1,a_1,\ldots ,a_{n-1}]+\sum_{i=1}^{n-1}(-1)^{i}y_1 a_i [1, h y_2,a_1,\ldots,\hat{a}_i,\ldots ,a_{n-1}]\\
     & \stackrel{\mathrm{\eqref{leib w 1}}}{=} y_1 h y_2 [1,a_1,\ldots ,a_{n-1}]+\sum_{i=1}^{n-1}(-1)^{i} y_1 a_i h [1, y_2,a_1,\ldots,\hat{a}_i,\ldots ,a_{n-1}]\\
     &+\sum_{i=1}^{n-1}(-1)^{i} y_1 a_i y_2 [1, h,a_1,\ldots,\hat{a}_i,\ldots ,a_{n-1}].
 \end{align*}

 Similarly,
 \begin{align*}y_2 [h y_1,a_1,\ldots,a_{n-1}] &\stackrel{\mathrm{\eqref{br as lcomb 1}}}{=}  y_2 h y_1 [1,a_1,\ldots ,a_{n-1}]+\sum_{i=1}^{n-1}(-1)^{i}y_2 a_i [1, h y_1,a_1,\ldots,\hat{a}_i,\ldots ,a_{n-1}]\\
 & \stackrel{\mathrm{\eqref{leib w 1}}}{=}  y_1 h y_2 [1,a_1,\ldots ,a_{n-1}]+\sum_{i=1}^{n-1}(-1)^{i} y_2 a_i h [1, y_1,a_1,\ldots,\hat{a}_i,\ldots ,a_{n-1}]\\ 
 &+\sum_{i=1}^{n-1}(-1)^{i} y_1 a_i y_2 [1, h,a_1,\ldots,\hat{a}_i,\ldots ,a_{n-1}].\end{align*}

Finaly,  applying  \eqref{br as lcomb 1} to $\sum_{i=1}^{n-1}(-1)^{i-1}h a_i[y_1,y_2,a_1,\ldots,\hat{a}_i,\ldots,a_{n-1}]$  we have 
  \begin{align*}\sum_{i=1}^{n-1}(-1)^{i-1}h a_i[y_1,y_2,a_1,\ldots,\hat{a}_i,\ldots,a_{n-1}] & \stackrel{\mathrm{\eqref{br as lcomb 1}}}{=} \sum_{i=1}^{n-1}(-1)^{i-1}h a_iy_1[1,y_2,a_1,\ldots,\hat{a}_i,\ldots,a_{n-1}]\\
& -\sum_{i=1}^{n-1}(-1)^{i-1}h a_iy_2[1,y_1,a_1,\ldots,\hat{a}_i,\ldots,a_{n-1}]\\
+\sum_{i=1}^{n-1}(-1)^{i-1}\sum_{j=1,(j\neq i)}^{n-1}(-1)^{j-1} & h a_i a_j[1,y_1,y_2,a_1,\ldots,\hat{a}_j,\ldots,\hat{a}_i,\ldots,a_{n-1}]=\end{align*}
(by \eqref{id: ha_ia_j[]} the last sum is equal to zero)
    $$=\sum_{i=1}^{n-1}(-1)^{i-1}h a_iy_1[1,y_2,a_1,\ldots,\hat{a}_i,\ldots,a_{n-1}]-\sum_{i=1}^{n-1}(-1)^{i-1}h a_iy_2[1,y_1,a_1,\ldots,\hat{a}_i,\ldots,a_{n-1}].$$

By substituting the previously derived expressions into the following sum, we obtain

  $$y_1 [h y_2,a_1,\ldots ,a_{n-1}]-y_2 [h y_1,a_1,\ldots,a_{n-1}]\\
    +\sum_{i=1}^{n-1}(-1)^{i-1}h a_i[y_1,y_2,a_1,\ldots,\hat{a}_i,\ldots,a_{n-1}]=0$$
in $(A,\cdot, [\cdot,\ldots,\cdot])$ for any $y_1,y_2,a_i\in A,$ and $1\leq i\leq n-1.$ Therefore, $(A,\cdot, [\cdot,\ldots,\cdot])$  is a strong  transposed Poisson $n$-Lie algebra.
\end{Proof}

In \cite{DzhConf}, it is noted that the algebra $(A, \cdot, W(\cdot, \ldots, \cdot))$ is a transposed Poisson $n$-Lie algebra. As a consequence of the theorem above, we establish that the transposed Poisson $n$-Lie algebra  $(A, \cdot, W(\cdot, \ldots, \cdot))$ is indeed a strong transposed Poisson $n$-Lie algebra.

\begin{corollary}\label{Dzh alg} Let
 $(A,\cdot, W(\cdot,\ldots,\cdot))$  be a commutative associative algebra $A$ with  $n-1$  commuting derivations $D_1,\ldots,D_{n-1}$, where the $n$-ary bracket is defined by
$$W(u_1,u_2,\ldots,u_n)=\begin{array}{|ccc|}
    u_1   &  \cdots &  u_n \\
  D_1 u_1   &  \cdots & D_1 u_n \\
    \cdots & \cdots & \cdots \\
      D_{n-1} u_1   &  \cdots & D_{n-1} u_n
\end{array},$$ 
for every $u_1,u_2,\ldots,u_n\in A.$
 Then $(A,\cdot, W(\cdot,\ldots,\cdot))$ is a strong transposed Poisson $n$-Lie algebra.
\end{corollary} 
\begin{Proof} Let $u_1,\ldots,u_{n}\in A.$ Then by definition of $W(\cdot,\ldots,\cdot)$ we have 
$$W(u_1,u_2,\ldots,u_n)=\sum_{i=1}^n(-1)^{i-1}u_i W(1, u_1,u_2,\ldots, \hat{u_i},\ldots,u_n).$$ 

Note that $$W(1,u_1 u_2,\ldots,u_n)=Jac(u_1 u_2,\ldots,u_n).$$ Hence, by \eqref{id: Jac leib}  we have 
$$W(1,u_1 u_2,\ldots,u_n)=u_1W(1, u_2,\ldots,u_n)+u_2W(1, u_1,\ldots,u_n) \; \; (n>2).$$  

By Theorem \ref{th: if bracket with 1 the Tr} the statement follows.
    
\end{Proof}

\section{\label{sec: simplicity TP}  Simplicity of transposed Poisson $n$-Lie algebras}

In this section, using the technique developed in \cite{Ouar2023}, we establish that a transposed Poisson $n$-Lie algebra is simple if and only if its associated $n$-Lie bracket $ [\cdot, \ldots, \cdot] $ is simple. The key distinction between this result for transposed Poisson $ 2 $-Lie algebras and transposed Poisson $n$-Lie algebras with $ n > 2 $ lies in the proof of the following proposition.

\begin{proposition}\label{id: I inside} The following identity holds true in a transposed Poisson $n$-Lie algebra $(A,\cdot, [\cdot,\ldots,\cdot]):$
    \begin{align*}
        \sum_{j=2}^n [[x_1,  \dots, x_{n}],y_2, \dots, y_j h,\dots, y_{n}]=
    \sum_{j=1, j\neq i}^n \sum_{i=1}^n (-1)^{i-1} [[x_i,y_2, \dots, y_{n}],x_1,  \dots, x_j h,\dots, \hat{x_{i}}, \ldots,x_n]. 
    \end{align*}  
    where $h,x_i,y_i\in  A$ and $i=1,\ldots,n.$
\end{proposition}
\begin{Proof}
   By the generalized Jacobi identity \eqref{id: jacobi} and the transposed Poisson identity \eqref{id: tr P}, we obtain 
    \begin{align*} 0&\quad  =n h([[x_1,  \dots, x_{n}],y_2, \dots, y_{n}]-\sum_{i=1}^n (-1)^{i-1} [[x_i,y_2, \dots, y_{n}],x_1, \dots, \hat{x_{i}}, \ldots,x_n])\\
    &\quad  = [x_1, \dots, x_{n}]h,y_2, \dots, y_{n}]+\sum_{j=2}^n [[x_1,  \dots, x_{n}],y_2, \dots, y_j h,\dots, y_{n}] \\  
    &\quad - \sum_{i=1}^n (-1)^{i-1} [[x_i,y_2, \dots, y_{n}]h,x_1, \dots, \hat{x_{i}}, \ldots,x_n] \\  
    &\quad -\sum_{j=1, j\neq i}^n \sum_{i=1}^n (-1)^{i-1} [[x_i,y_2, \dots, y_{n}],x_1,  \dots, x_j h,\dots, \hat{x_{i}}, \ldots,x_n].
\end{align*}  
By identity \eqref{id: [[]h]}, the first and third terms gives zero, and we have 
     $$\sum_{j=2}^n [[x_1,  \dots, x_{n}],y_2, \dots, y_j h,\dots, y_{n}]=\sum_{j=1, j\neq i}^n \sum_{i=1}^n (-1)^{i-1} [[x_i,y_2, \dots, y_{n}],x_1,  \dots, x_j h,\dots, \hat{x_{i}}, \ldots,x_n].$$ 
     This completes the proof.

\end{Proof}

Let $(A,\cdot, [\cdot,\ldots,\cdot])$ be a transposed Poisson $n$-Lie algebra. If the $n$-Lie algebra $(A, [\cdot,\ldots,\cdot])$ satisfies $[A,\ldots,A] \neq A$, then by identity \eqref{id: tr P}, the subspace $[A,\ldots,A]$ is an ideal of the commutative associative algebra $(A,\cdot)$. Moreover, it is also an ideal of the transposed Poisson $n$-Lie algebra $(A,\cdot, [\cdot,\ldots,\cdot])$. Therefore, in any simple transposed Poisson $n$-Lie algebra $(A,\cdot, [\cdot,\ldots,\cdot])$, the associated $n$-Lie algebra $(A,[\cdot,\ldots,\cdot])$ must satisfy $[A,\ldots,A] = A$. 

Recall that a quasi-ideal of a transposed Poisson $n$-Lie algebra $(A,\cdot, [\cdot,\ldots,\cdot])$ is a proper subspace $I \subset A$ such that  
\[
[I,A,\ldots,A] \subset I \quad \text{and} \quad [AI,A,\ldots,A] \subset I.
\]

The following two lemmas are generalizations of Lemma 9 and Lemma 10 in \cite{Ouar2023}. 

\begin{lemma}\label{lem: no quasi-id} A simple transposed Poisson $n$-Lie algebra $(A,\cdot, [\cdot,\ldots,\cdot])$ does not contain quasi-ideals.
\end{lemma}
\begin{Proof}
    Suppose that $I$ is a quasi-ideal of $A.$ Let $I'$ be a maximal subspace of $A$ such that $[I',A,\ldots,A]\subset I.$ Since $(A,\cdot, [\cdot,\ldots,\cdot])$ is simple, we have $[A,A,\ldots,A]=A,$ so $I'\neq A.$ Moreover, by maximality of $I'$ and since $I$ is quasi-ideal of $A,$ we have $AI\subset I'.$ Now, by identity \eqref{id: h[)} we have
    $$ I'\, [a_1,\ldots,a_n]\subset \sum_{i=1}^{n} a_i \, [a_1,\ldots,a_{i-1}, I',a_{i+1},\ldots, a_n]\subset I',$$
     for every $a_1,\ldots,a_n\in A.$ Since $[A,A,\ldots,A]=A,$ we have $AI'\subset I',$ and  $[I',A,\ldots,A]\subset I\subset I'.$
     Thus, $I'$ is an ideal of $A$ which contradicts the simplicity of $(A,\cdot, [\cdot,\ldots,\cdot]).$
\end{Proof}

\begin{lemma}\label{lem: id is quasi-id}
    Let $(A,\cdot, [\cdot,\ldots,\cdot])$ be a transposed Poisson $n$-Lie algebra and $[A,\ldots,A]=A.$ Then any ideal in the $n$-Lie algebra $(A,[\cdot,\ldots,\cdot])$ is a quasi-ideal.
\end{lemma}
\begin{Proof}
    Suppose that $I$ is an ideal of the $n$-Lie algebra $(A,[\cdot,\ldots,\cdot]).$ By Proposition  \ref{id: I inside} we have
     $$[[x_1,  \dots, x_{n}],I h, y_3,\dots, y_{n}]\subset \sum_{j=3}^n [[x_1,  \dots, x_{n}],I,y_3, \dots, y_j h,\dots, y_{n}]+$$$$\sum_{j=1, j\neq i}^n \sum_{i=1}^n (-1)^{i-1} [[x_i,I, \dots, y_{n}],x_1,  \dots, x_j h,\dots, \hat{x_{i}}, \ldots,x_n]\subset I,$$ for every $h,x_i,y_i\in A.$
    Hence, $I$ is a quasi-ideal.
\end{Proof}

By applying the two lemmas above, we can derive the following theorem.

\begin{theorem}
    Let $(A,\cdot, [\cdot,\ldots,\cdot])$ be a simple transposed Poisson $n$-Lie algebra. Then $(A, [\cdot,\ldots,\cdot])$ is simple.
\end{theorem}
\begin{Proof} Suppose that $(A, [\cdot,\ldots,\cdot])$ is not simple, then any ideal is a quasi-ideal by Lemma \ref{lem: id is quasi-id}, but this contradicts Lemma \ref{lem: no quasi-id}, because $A$ is simple. Therefore, the $n$-Lie algebra $(A, [\cdot,\ldots,\cdot])$ must be simple.
\end{Proof}

Therefore, we have the main result of this section:

\begin{corollary}\label{siffs}
     Let $(A,\cdot, [\cdot,\ldots,\cdot])$ be a transposed Poisson $n$-Lie algebra. Then $(A, \cdot, [\cdot,\ldots,\cdot])$ is simple if and only if $(A,  [\cdot,\ldots,\cdot])$ is simple.
\end{corollary}

\begin{remark}
   In \cite{DzhJac}, it was proved that $(A, W(\cdot,\ldots,\cdot))$ is a simple $n$-Lie algebra. Therefore, by Corollary \ref{Dzh alg} and Corollary \ref{siffs} we have $(A,\cdot, W(\cdot,\ldots,\cdot))$ is a simple transposed Poisson $n$-Lie algebra. As noted in the introduction the algebra $(A,\cdot, W(\cdot,\ldots,\cdot))$ is a generalized Nambu-Poisson algebra that appears in the classification presented in \cite{CanKac2016}.
\end{remark}

\section{A strong transposed Poisson $3$-Lie algebras}

In this section, we use linear algebra to demonstrate that the strong condition does not hold in the free transposed Poisson 3-Lie algebra.

Recall that the Leibniz identity \eqref{id: tr P} for $n=3$ has the following form
$$T(h,a_1,a_2,a_3)= 3\, h [a_1,a_2,a_3]- [h a_1,a_2, a_3]- [a_1,h a_2, a_3]- [a_1,a_2,h a_3]=0,$$
and the strong condition:
$$S(h,y_1,y_2,x_1,x_2)=y_1 [h y_2,x_1,x_{2}]-y_2 [h y_1,x_1,x_{2}]
    +h x_1[y_1,y_2,x_2]-h x_2[y_1,y_2,x_1]=0.$$

Now, let us consider the free transposed Poisson $3$-Lie algebra $\mathcal{A}(X)$ generated by  a countable set $X = \{a_1, a_2, \ldots\}$.  We will prove that this algebra does not satisfy the strong condition.

\begin{proposition}
    Let $\mathcal{A}(X)$
be a free transposed Poisson $3$-Lie algebra. Then $\mathcal{A}(X)$ does not satisfy the identity $S(h,y_1,y_2,x_1,x_2)=0$.
\end{proposition}

\begin{Proof}

Since $T(h,a_1,a_2,a_3)$ is skew-symmetric in its last three arguments, it follows that $$T(h,a_1,a_2,a_3)=sign(\sigma) T(h,a_{\sigma(1)},a_{\sigma(2)},a_{\sigma(3)})$$ for $\sigma\in S_3.$
From this, we obtain three consequences of the Leibniz identity $T(h,x,y,z)=0$ at degree 5, by multiplying the identity and arguments by the associative product $``\cdot":$
\begin{equation}\label{cons of T-3}
    T(a_1, a_2, a_3, a_4)a_5,\,T(a_1 a_5, a_2, a_3, a_4),\,T(a_1 , a_2 a_5, a_3, a_4).\end{equation}

We apply all permutations $\sigma \in S_5$ to the arguments $a_1,\ldots,a_5$  of the polynomials \eqref{cons of T-3}.
Since the product $ \cdot$  is symmetric and the bracket \( [\cdot,\cdot,\cdot] \) is skew-symmetric, we obtain a set of 65 monomials of degree 5. We present these monomials  in the following order:
\begin{equation}\label{monomials of d 6}\begin{array}{c}
   \{a_{\sigma(1)}a_{\sigma(2)}[a_{\sigma(3)},a_{\sigma(4)},a_{\sigma(5)}]|\; \sigma(1)<\sigma(2),\; \sigma(3)<\sigma(4)<\sigma(5), \sigma \in S_5\},\\
\{a_{\sigma(1)}[a_{\sigma(2)} a_{\sigma(3)},a_{\sigma(4)},a_{\sigma(5)}]|\; \sigma(2)< \sigma(3),\sigma(4)<\sigma(5), \sigma \in S_5\},\\
\{[a_{\sigma(1)}a_{\sigma(2)},a_{\sigma(3)}a_{\sigma(4)},a_{\sigma(5)}]|\; \sigma(2)>\sigma(1)< \sigma(3)<\sigma(4), \sigma \in S_5\},\\
\{[a_{\sigma(1)} a_{\sigma(2)} a_{\sigma(3)},a_{\sigma(4)},a_{\sigma(5)}]|\; \sigma(1)<\sigma(2)<\sigma(3),\; \sigma(4)<\sigma(5), \sigma \in S_5\}.  
\end{array}
\end{equation}

Then, each polynomial in \eqref{cons of T-3} is a linear combination of monomials \eqref{monomials of d 6}. The skew-symmetry of $T(h,x,y,z)$ in arguments $x,y,z,$ allows us to obtain the row space  $120 \times 65$ matrix $C$ whose columns are labelled by the monomials \eqref{monomials of d 6}. Computing its rank, we find that $rank(C)=46.$

Now, suppose that a free transposed Poisson 3-Lie algebra $\mathcal{A}(X)$ satisfies the strong condition:
$S(h,y_1,y_2,x_1,x_2)=0.$ It is easy to see that the polynomial $S(a_1,a_2,a_3,a_4,a_5)$ cannot be written as a linear combination of the generalized Jacobi polynomials.
Then, $S(a_1,a_2,a_3,a_4,a_5)$ must be represented as a~linear combination of polynomials \eqref{cons of T-3}. 

However, if we add an additional row corresponding to $S(a_1,a_2,a_3,a_4,a_5)$
 to our matrix $C$, we obtain a new matrix $C'$ of size $121 \times 65.$ Computing its rank, we find $rank(C')=47.$ Thus, we have a contradiction which proves that a free transposed Poisson 3-Lie algebra is not strong. 
\end{Proof}


\begin{center} ACKNOWLEDGMENTS   \end{center}
This research was funded by the Science Committee of the Ministry of Science and Higher Education of the Republic of Kazakhstan (Grant No. AP22683764).

\textbf{There is no conflict of interest}

\end{document}